\documentclass[11pt]{article}

\usepackage{latexsym}
\usepackage{amssymb}

\renewcommand{\d}{\mathrm{d}}

\newcommand{\koniec}{\begin{flushright}  $\Box $ \end{flushright}}
\def\be{\begin{equation}}
\def\ee{\end{equation}}

\def\t{\tilde}

\def\Om{\Omega}

\def\om{\omega}

\def\p{\partial}

\newtheorem{prop}{Proposition}  
\newtheorem{lemma}[prop]{Lemma}
\begin{document}
\title{Scalar--Flat Lorentzian Einstein--Weyl Spaces}
\author{David M.J.Calderbank\\
Department of Mathematics and Statistics,\\
 University of Edinburgh,\\ 
King's Buildings, Mayfield Road, Edinburgh EH9 3JZ, Scotland UK
\and
Maciej  Dunajski\\ 
Mathematical Institute,\\
University of Oxford, 24-29 St Giles, Oxford OX1 3LB, UK.}  
\date{} 
\maketitle
{\abstract{We find all three-dimensional Einstein--Weyl spaces 
with the vanishing scalar curvature.}} 
\vskip5pt
\noindent
Three-dimensional Lorentzian Einstein metrics have constant curvature.
In order to construct  
nontrivial gravitational models in three dimensions, one should therefore 
look at conformal geometries involving  a non-metric connection.
Einstein--Weyl geometries appear quite naturally in this context.

Let ${\cal W}$ be a
three-dimensional Weyl space i.e. a manifold with 
with a torsion-free connection $D$ and a
conformal metric $[h]$ such that 
null geodesics of $[h]$ are also geodesics for $D$.  This condition is
equivalent to 
$
D_ih_{jk}=\om_ih_{jk}
$
for some one form $\om$. Here $h_{jk}$ is a representative metric in
the conformal class, and the indices 
$i, j, k, ...$ go from $1$ to $3$.
If we change this representative by
$h\longrightarrow \phi^2 h$, then $\om\longrightarrow
\om+2\d\ln{\phi}$.  
The one-form $\om$ `measures' the difference between $D$ and
the Levi-Civita connection $\nabla$ of $h$.
A tensor object $T$ which transforms as $
T\longrightarrow \phi^m T$ when
$ h_{ij}\longrightarrow \phi^2 h_{ij}$
is said to be conformally invariant of weight $m$.
The formula for a
covariant weighted derivative of a one-form of weight $m$ is
\be
\label{wcowder}
{D}_iV_j= \nabla_iV_j+\frac{1}{2}
((1-m)\om_iV_j+\om_jV_i-h_{ij}\om_kV^k).
\ee
The Ricci tensor $W_{ij}$  
is related to the Ricci tensor
$R_{ij}$ and of $\nabla$ by
\[
W_{ij}=R_{ij}+\nabla_i\om_j-\frac{1}{2}\nabla_j\om_i
+\frac{1}{4}\om_i\om_j+h_{ij}\Big(-\frac{1}{4}\om_k\om^k+\frac{1}{2}
\nabla_k\om^k\Big).
\]
The conformally invariant Einstein--Weyl (EW) condition on
$({\cal W}, h, \om)$ is
$
W_{(ij)}=W h_{ij}/3,
$
or in  terms of the Riemannian data: 
\be
\label{ew2}
\chi_{ij}:=R_{ij}+\frac{1}{2}\nabla_{(i}\om_{j)}+\frac{1}{4}\om_i\om_j
-\frac{1}{3}\Big(r+\frac{1}{2}\nabla^k\om_k+\frac{1}{4}\om^k\om_k\Big)h_{ij}=0,
\ee
where $\chi_{ij}$ is the trace-free part
of the Ricci tensor of the Weyl connection, and $r=h^{ij}R_{ij}$.
Weyl spaces which satisfy (\ref{ew2}) will be called Einstein--Weyl
spaces.
The EW equations can be regarded as an integrable system. This is
because
both the twistor theory \cite{Hi82} and the Lax representation
\cite{DMT00} exist. One should therefore be able to construct 
large families of explicit solutions. 

In this paper we shall find explicitly all EW spaces with vanishing
scalar curvature $W=h^{ij}W_{ij}$. 
We shall first establish the following result:
\begin{lemma}
\label{nullF}
If the scalar curvature of the Weyl connection vanishes on and EW
space $({\cal W}, [h], D)$, then
the Faraday two form $F_{ij}:=\nabla_{[i}\om_{j]}$ is null.
\end{lemma}
{\bf Proof.}
The Bianchi identities for the curvature of the Weyl connection
written in terms of the Levi-Civita connection and $\om$ are \cite{PT93}
\be
\label{bianchi}
\nabla^iF_{ij}+\frac{1}{2}\om^iF_{ij}+\frac{1}{3}(\nabla_j W
+\om_jW)=0.
\ee
Assume $W=0$ (this is a well defined  condition as $W$ is
conformally invariant of weight $-2$).
Contracting (\ref{bianchi}) with $\nabla^j$, and using
$\nabla^j\nabla^iF_{ij}=0$ yields
\[
0=(\nabla^j\om^i)F_{ij}+\om^i\nabla^jF_{ij}=F^{ij}F_{ij}-
\frac{1}{2}\om^i\om^jF_{ij},
\]
so $F$ is null.\koniec
We conclude that there are no non-trivial scalar-flat EW spaces  
in the Euclidean signature \cite{Ca98}.
Non-trivial solutions can be found in the indefinite signature:
\begin{prop}
Let $(h, \om)$ be an Einstein--Weyl structure with vanishing scalar
curvature. Then either  $(h, \om)$  is flat, or the signature is $(++-)$
and there exist local coordinates  $x^i=(y, x, t)$ 
such that $\om=y\d t$, and $h$ 
is given by one of two solutions:
\be
\label{case1}
h_1=\d y^2+2\d x\d t+ \Big(x[R(t)-\frac{y}{2}]+\frac{1}{48}y^4+\frac{1}{12}R(t)y^3
+S(t)y\Big)\d t^2,
\ee
\be
\label{case2}
h_2=\d y^2+2\d x\d t-\frac{4x}{y}\d y\d t +\Big(\frac{x^2}{y^2}
+\frac{xy}{2}+\frac{1}{8}y^4+R(t)y^2
+S(t)y\Big)\d t^2,
\ee
where $R(t)$ and $S(t)$ are arbitrary functions with 
continuous second derivatives.
\end{prop}
{\bf Proof.}
Lemma \ref{nullF} implies  that $F_{ij}=\nabla_{[i}\om_{j]}$ is a closed
null two-form. The conformal freedom together with the Darboux theorem imply
the existence of coordinates such that $\om_i =y\nabla_i t$. 
Therefore $\om\wedge\d \om=0$, and the nullity of $F$ gives 
$\ast F\wedge\om=0$. We can rewrite
Bianchi identity (\ref{bianchi}) as
\[
2(\d \ast F)+\om\wedge\ast F=0,
\]
and deduce $\d \ast F=0$.
Therefore
$\nabla^{i}F_{ij}=0$, 
and ${\varepsilon_{i}}^{jk} F_{jk}=f(t)\nabla_i t$. 
Redefining $y, t$ we can set $f(t)=1$.
The most general metric consistent with $\d t=\ast\d y\wedge\d t$ is
\[
h=\d y^2+2(\hat{E}\d s+\hat{F}\d y+\hat{G}\d t)\d t, 
\]
where $\hat{E}, \hat{F}$, and $\hat{G}$ 
are functions of $(s, y, t)$. Put $\hat{E}=\p x/\p s$
and define $G=(\hat{G}-x_t)/2, F=\hat{F}-x_y$, so that
\be
\label{EW1}
h=\d y^2+2\d x\d t+2F\d y\d t+G\d t^2,\qquad \om=y\d t.
\ee
The freedom $x\rightarrow x+P(y, t)$ implies that 
$F(x, y, t)$, and $G(x, y, t)$ are defined up to addition of
derivatives of $P(y, t)$. Furthermore the conformal scale is only
fixed up to arbitrary functions of $t$, $h\mapsto\tilde{h}=\Omega h$. 
This leads to  to the redefinitions $(x,y,t) \rightarrow
(\tilde{x},\tilde{y},\tilde{t})$, given by
\be\label{conftrans}
\t{t}
=T(t),\qquad
\t{y}=\frac{y}{T_t}-2\frac{T_{tt}}{T_t},\qquad\t{x}=\frac{x}{{T_t}^3}+P(y,t),
\qquad \Om=T_t^2.
\ee
These transformations will latter 
be used to simplify  $F$ and $G$. 
Now impose the EW equations:
Equation $\chi_{12}=0$ implies $F_{xx}=0$, and we can
choose $P(y, t)$ such that $F=xl(y,t)$. Now
$\chi_{12}=\chi_{22}=0$. Equations $\chi_{11}=0, \chi_{23}=0$ are equivalent
and imply $a_y=G_{xx}$. Take $G=x^2l_y/2+xm(y, t)+n(y, t)$. 
The vanishing of the scalar curvature 
\[
W=r+2\nabla^k\om_k-\frac{1}{2}\om^k\om_k=(3l^2-6l_y)/2
\]
gives $l^2=2l_y$, therefore
$l(y, t)=0$ (case {\bf 1}), or $l(y, t)=-2/(y+c(t))$ (case {\bf 2}). 
\begin{itemize}
\item
In case {\bf 1} $\chi_{31}=0$ implies $m(y, t)=R(t)-y/2$.
Finally $\chi_{33}=0$ yields $n(y, t)=y^4/48+R(t)y^3/12
+S(t)y+Z(t))$. Redefining $x$, we  set $Z(t)=0$, and the metric is
given by (\ref{case1}).
\item In case {\bf 2}  the conformal freedom 
(\ref{conftrans}) can be used to eliminate $c(t)$. 
This is achieved by setting $T_{tt}=c(t)$, and redefining $m(y, t)$, 
and $n(y, t)$.
Now  $\chi_{31}=0$ implies $m(y, t)=y/2 +P(t)$, and
$\chi_{33}=0$ gives $n(y, t)=y^4/8+y^3P(t)/4 +R(t)y^2+S(t)y$.
In fact we can get rid of $P(t)$: Replace $x$ by $x-y^2w(t)/2$, and
redefine $R(t)$ to obtain (\ref{case2}).
\end{itemize} 
\koniec
The next proposition shows that there doesn't exist a combination 
of coordinate and conformal transformations which maps $(h_1, \om)$
to $(h_2, \om)$. Cases {\bf 1} and {\bf 2} are essentially distinct
and can be invariantly characterized:
\begin{prop}
Let $(h, \om)$ be a non-flat EW structure  with a vanishing
scalar curvature, and let $F$ be a corresponding Faraday two form.
If $\ast F$ is parallel with respect to a weighted Weyl connection, 
then $(h, \om)$ is locally given by
{\em(}\ref{case1}{\em)}. Otherwise it is locally given by {\em (}\ref{case2}{\em )}.
\end{prop}
{\bf Proof.}
In both cases {\bf 1} and {\bf 2}
$\ast F={\varepsilon_i}^{jk}F_{jk}\d x^i =\d t$. First
notice that vanishing of $D(f\d t)$ for some $f$ is invariant
condition; in a conformal scale defined by $f$ the one-form $\ast F$
is parallel. Here we treat $f\d t$ as a weighted object.
In case {\bf 1} we find that with $m=3/2$ we have $D(f\d t)=0$ if $f=f(t)$, and
$f_t+fR/2=0$.
In case {\bf 2}  $D(f\d t)$ doesn't vanish for any $f$. For example
$(D(f\d t))_{13}=-f/y$.
\koniec
Using the formula for a weighted derivative of a vector of weight $m$
\[
D_iV^j=\nabla_iV^j-\frac{1}{2}\delta_i^j\om_kV^k
-\frac{m+1}{2}\om_iV^j+\frac{1}{2}\om^jV_i
\]
we deduce that the EW structure (\ref{case1}) (case {\bf 1})
admits a covariantly constant vector
\[
V=\exp\Big(-\frac{1}{2}\int^tR(t)\d t\Big)\frac{\p}{\p x}
\]
with weight $-1/2$.  Therefore 
$(h_1, \om)$ belongs to the conformal 
class of the dKP Einstein--Weyl spaces \cite{DMT00}.

It is natural to ask if other special classes of Einstein--Weyl spaces
include Lorentzian scalar-flat examples. Along these lines, we have
the following observation.
\begin{prop} Let $R$ be an arbitrary function of one variable. Then
the Weyl structure
\be
h= 4\frac{(z+R(v))^2}{(1+vw)^2}\d v \d w +\d z^2, \qquad
\om = \frac 4{z+R(v)} \d z
\ee
is scalar-flat, and $u=2\log \bigl(2(z+R(v))/(1+vw)\bigr)$ is a
solution of the Lorentzian $SU(\infty)$ Toda equation $4u_{vw}+(e^u)_{zz}=0$.
\end{prop}
{\bf Proof.}
These are straightforward verifications, using the fact that the
scalar curvature of the Weyl structure ($h=e^u\d v\d w+\d z^2$,
$\om=2u_z\,\d z$) is $\frac12 u_{zz}+\frac14 u_z^2$:
these spaces are Lorentzian analogues of the hyperCR--Toda
Einstein--Weyl spaces~\cite{CT01}, see also \cite{Wa90}.
\koniec


\end{document}